\documentclass[12pt]{article}
\usepackage{setspace}
\usepackage{textcomp}
\usepackage{amsmath, amsthm, amssymb}
\usepackage[left=1.5in,top=.9in, 
right=1in,nohead]{geometry}
\usepackage{latexsym}
\usepackage{mathrsfs}
\usepackage{rawfonts}
\usepackage[dvips]{graphicx}
\usepackage[all]{xy}
\usepackage[sc,osf]{mathpazo}  

\def\bct{\begin{center}}
\def\ect{\end{center}}
\def\beg{\begin}

\def\<{\langle}
\def\>{\rangle}

\def\mbb{\mathbb}
\def\mbbr{\mathbb R}
\def\mbbz{\mathbb Z}

\def\mc{\mathcal}

\def\ni{\noindent}
\def\ol{\overline}

\def\tn{\textnormal}
\def\wt{\widetilde}
\def\ti{\tilde}

\newtheorem{thm}{Theorem}[section]

\newtheorem{prop}[thm]{Proposition}
\newtheorem{cor}[thm]{Corollary}

\title{
Free immersions and panelled web 4-manifolds} \author{Mustafa Kalafat}

\begin{document}
\maketitle
\begin{abstract} We show that if a compact, oriented 4-manifold admits a 
coassociative($*\phi_0$)-free immersion into $\mbb R^7$ then its 
Euler characteristic $\chi_M$ and signature 
$\tau_M$ vanish. Moreover, in the spin case the Gauss map 
is contractible, so that the immersed manifold is parallelizable. 
The proof makes use of homotopy theory in particular obstruction theory.  
As a further application we prove a non-existence result for some infinite families of 4-manifolds that can not be addressed previously. 
We give concrete examples of parallelizable 4-manifolds with 
complicated non-simply-connected topology.

\end{abstract}
\section{Introduction}

Let $(M,g)$ be a Riemannian $n$-manifold. If we take a point $p$ in $M$, a $k$-dimensional 
vector subspace $V<T_pM$ equipped with an orientation is called an {\em oriented tangent $k$-plane} of $M$. In this case the restricted metric $g|_V$ and the orientation gives a natural k-form, the volume form $\tn{vol}_V$ on $V$. A $k$-form $\varphi$ on $M$ is called a {\em calibration} if it is closed and $\varphi|_V \leq \tn{vol}_V$ for any oriented $k$-plane $V$. 
In general $\varphi|_V=\alpha\cdot \tn{vol}_V$ for some $\alpha\in\mbb R$ since they are both top forms on the vector space $V$. The calibration condition is equivalent to  $\alpha\leq 1$. Under these assumptions if $N$ is an oriented $k$-dimensional submanifold of $M$, then the tangent spaces of $N$ are automatically oriented tangent $k$-planes and we say that $N$ is a {\em calibrated submanifold} or $\varphi$-{\em submanifold} of $M$ if $N$ has maximal tangent spaces i.e. $\varphi |_{T_pN}=\tn{vol}_{T_pN}$ for all $p\in N$. The function $\alpha\equiv 1$ constant on $N$ in this case. Calibrated manifolds are introduced in \cite{harveylawson}, for a survey see also \cite{joyce}. 
On the Euclidean space $\mbb R^7$, consider the {\em coassociative form} which is the 4-form
$$*\phi_0=dx^{1234}-dx^{12-34}\wedge dx^{67}-dx^{13-42}\wedge dx^{75}
                   -dx^{14-23}\wedge dx^{56}.$$
This is actually an example of a calibration. 
The subgroup of $GL(7,\mbb R)$ that leaves $*\phi_0$ invariant is the compact $14$-dimensional Lie group $G_2$. If $U$ is the $4$-plane in $\mbb R^7$ with last 
 three coordinates vanishing, then $*\phi_0|_U=dx^{1234}$ which is equal to $\tn{vol}_V$ with suitable orientation on $U$ hence $U$ is a maximal plane called a {\em coassociative $4$-plane}. One can show that \cite{harveylawson} the subgroup of $G_2$ preserving $U$ is $SO(4)$, and $G_2$ acts transitively on (oriented) coassociative $4$-planes. So that the set COASS of coassociative planes in $\mbb R^7$ is isomorphic to $G_2/SO(4)$ and has dimension $8$. The Grassmannian $G_4^+\mbb R^7$ of all oriented 4-planes in $\mbb R^7$ is of dimension $12$ so COASS is a codimension 4 subspace of it. We also have that $\phi_0|_V=0$ for every coassociative 4-plane $V$ since $\phi_0|_U=0$ by definition, the action of $G_2$ is transitive on the coassociative 4-planes and $\phi_0$ is $G_2$ invariant. Conversely if 
$\phi_0|_V=0$ for a 4-plane then there is a unique orientation which makes $V$ a coassociative 4-plane.

\vspace{.05in}

We are actually interested in the following type of submanifolds. Following \cite{harveylawsonpotential} if $(M,\phi)$ is a calibrated manifold, a submanifold $N$ is called {\em $\phi$-free} if there are no $\phi$-planes tangent to $M$. 
These submanifolds have strictly $\phi$-convex neighborhoods each of which admits deformation retraction onto $N$. These are generalizations of totally real submanifolds in K\"ahler manifolds to calibrated manifolds, and strictly $\phi$-convex manifolds are the generalization of Stein manifolds. They have nice topological structures. 
In order to understand $\phi$-free submanifolds, some special information is needed about the related Grassmann manifolds. 
We will be using the 
results of \cite{atg2} on the topology of Grassmannians.
As an application 
we will give an answer to  
coassociative-free embedding problems for some infinite families of 4-manifolds. 
There is very few results on the coassociative-free embeddings of 4-manifolds. Only notable obstruction is the Euler characteristic of 
\. I. \" Unal in \cite{ibram}. One can for example conclude that coassociative-free embeddings of the 4-manifolds $S^4$, $\Sigma_g\times \Sigma_h$ for genus $g,h\neq 1$ are violated since they have non-zero 
Euler characteristic. On the other hand, the techniques 
in the literature can not answer this question for $\chi=0$ case. 
Due to our main result in this paper, we now able find better obstructions as follows.

\vspace{.08in}

\noindent {\bf Theorem \ref{vanishingtheorem}.}(Vanishing).
{\em  If $M^4$ is closed and $i: M\to \mbb R^7$ is a
coassociative($*\phi_0$)-free immersion, then the Euler characteristic $\chi_M$ and 
the signature $\tau_M$ vanishes.}

\vspace{.07in} 

\ni If we combine it with the results in \cite{ibram3} on the converse, we obtain the following. 

\vspace{.05in}

\begin{cor}
A closed $4$-manifold $M$ admits a
coassociative($*\phi_0$)-free immersion or embedding into $\mbb R^7$ 
if and only if its Euler characteristic $\chi_M$ and 
the signature $\tau_M$ vanishes.\end{cor}

\ni One can also change the target space from the Euclidean space to any manifold with $G_2$ structure which is flat in a neighborhood of a point.  

\vspace{.05in}

In the spin case we are now able to give a better obstruction as follows. 
We will make use of the fact that in dimension four, 
parallelizability is characterized through the following complete obstructions.

\begin{thm}[\cite{hirzebruchhopf,masseyparallelizability}]\label{thmparallelizabilityHHM}
A smooth $4$-manifold is parallelizable iff  $w_{12}=e=p_1=0$. 
\end{thm}
\ni See also \cite{emeryparallelizability}. So this implies that, 
further in the spin case the $4$-manifold has to be parallelizable. 
We alternatively prove this fact going through the analysis of the Gauss map,
and show that it is trivializable as well.
\vspace{.08in}

\noindent {\bf Theorem \ref{contractible}.}
{\em If $M^4$ is a compact, oriented, spin manifold 
and $i: M\to \mbb R^7$ is a
coassociative($*\phi_0$)-free immersion, then its image $g(M)$ under the Gauss map $g: M\to G_4^+\mbb R^7$ is contractible, so that $M$ is parallelizable. 
}


\vspace{.05in}

Our Proposition \ref{image} is 
crucially used in the proof of the main result of the paper \cite{ibram3}. 
So our paper fills a gap in the literature in this perspective. 
As another application we illustrate our results through a series of examples in 
section \S\ref{secexamples}. 
We do a similar computation from this series in section \S\ref{seccayleyfree} 
for the Cayley-free case. In \cite{ibram4} using h-principle techniques, an 
if and only if theorem is proved for this case as well. 
In section \S\ref{embeddings} we prove our main result, 
in section \S\ref{hurewiczinjectivity}  we prove an injectivity lemma and applications 
on the Gauss map. 

\vspace{.05in}

{\bf Acknowledgements.} We thank B. Lawson 
and S. Akbulut for their 
 suggestions. 
Many thanks to J. Morgan for very useful remarks.     
Thanks to T. \" Onder for some referencing. Thanks to the anonymous referee for useful remarks.  This work is 
partially supported by T\"ubitak (Turkish science and research council) grant {$\sharp$}114F320.


\section{Coassociative-free Immersions}\label{embeddings}

In this section we will present applications to the immersion theory. 
We will be using results on the topology of the oriented Grassmann space $G_3^+\mbb R^7$. 
We denote the canonical(tautological) vector bundle and its orthogonal complement by $E=E_3^7$ and $F=F_3^7$ on this space. 
We will often be using the following result. 

\beg{thm}[\cite{shizhou}]\label{shizhouthm} We have the following characteristic class relations for the bundles over the Grassmannian $G_3^+\mbb R^7$,\\

(a) $p_1E=-p_1F$, $p_1^2E=e^2F$.\\

(b) $p_1E[\mbb{CP}_2]=p_1E[\overline{\mbb{CP}}_2]=eF[\mbb{CP}_2]=-eF[\overline{\mbb{CP}}_2]=1$.\\

(c) ${1\over 2}(p_1E \pm eF)$ are generators in $H^4(G_3^+\mbb R^7;\mbbz)$. Their Poincar\'e duals are\\

$[ASS]$ and $[\wt{ASS}]$ respectively.\\

(d) ${1\over 2}(p_1E eF \pm e^2F)$ are generators in $H^8(G_3^+\mbb R^7;\mbbz)$. Their 
Poincar\'e duals are\\

$[\mbb{CP}]$ and $[\overline{\mbb{CP}}_2]$ respectively.
\end{thm}

This is a combination of the results 
in Section 7 of the resource. 
Notice that, to be able to say that the generators stated in part (c) are the sole 
generators, one needs to know that there is no torsion as explained in \cite{atg2}. 
Also note that we can realize the embeddings of $\mbb{CP}_2$ and $\overline{\mbb{CP}}_2$ through 
the inclusion of $G_2^+\mbb R^6$ in $G_3^+\mbb R^7$ since the oriented Grassmannian is a double cover of the Grassmannian, we can include the projective space with both orientations. 
These two projective spaces are the generators of $H_4(G_3^+\mbb R^7;\mbb R)$ as explained in the resource. 

\vspace{.05in}
 
Let $i:M\to \mbb R^7$ be an 
immersion of a 4-manifold into the Euclidean space
\footnote{See \cite{fang} for a discussion on the topological embedding problem.} 
. In this case we have the associated Gauss map $g:M\to G_4^+\mbb R^7$ from the manifold to the oriented Grassmannian. It sends a point to the 4-dimensional subspace of the Euclidean space which is parallel to the tangent space at that point. 
In general no such map exist due to lack of translation unless e.g. the space is flat and simply connected or alternatively parallelizable. 
Composing this map with the orthogonal complement map $*:G_4^+\mbb R^7\to G_3^+\mbb R^7$ we get the map 
$$\tilde{g} : M^4\longrightarrow G_4^+\mbb R^7 \longrightarrow G_3^+\mbb R^7 ~~~\tn{where}~~~ \tilde{g}=* \circ g$$ 
which is appropriate setting for us to apply the results above. 
Keep in mind that the orthogonal complement map the pullbacks the bundles as  
$*^*E_3^7=F_4^7$ and $*^*F_3^7=E_4^7$.  
We start with our first lemma.

\vspace{1cm}

\begin{figure}[h] \bct \label{cointersectionnumber}
\includegraphics[width=.8\textwidth]{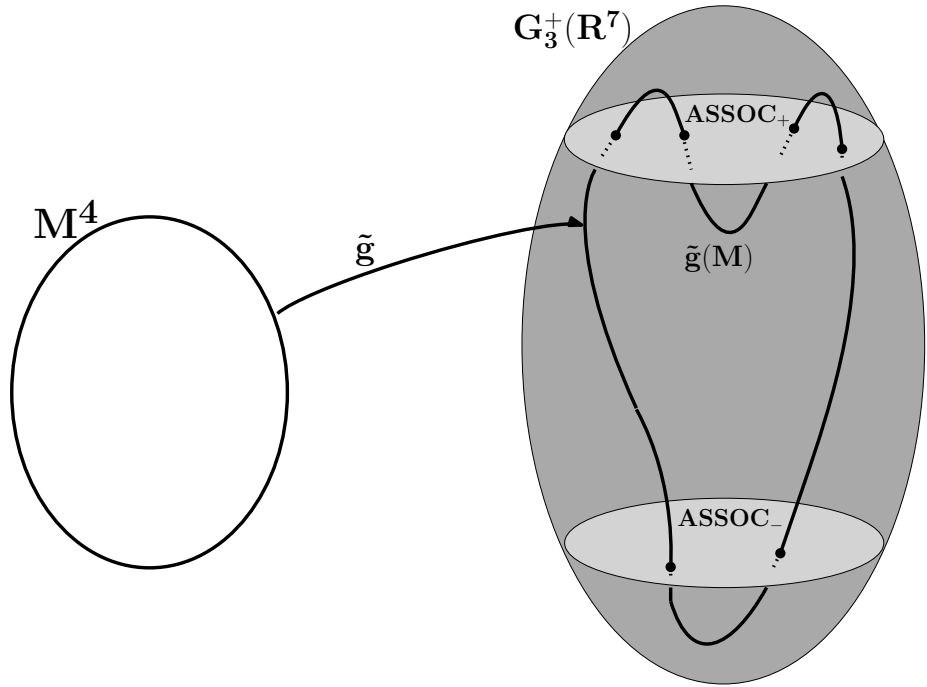}\\ \ect
  \caption{\, The Gauss map and intersections\,  } 
\end{figure}

\beg{lem}\label{c1} 
The image $\tilde g_*[M]=c\, [\mbb{CP}_2]+(c-\chi_M)  [\overline{\mbb{CP}}_2]
\in H_4(G_3^+\mbb R^7;\mbb R)$ for some $c\in\mbb R$.

\end{lem}
\beg{proof} Suppose $\tilde g_*[M]=c [\mbb{CP}_2]+d [\overline{\mbb{CP}}_2]$ for some  $c,d \in \mbb R$. 
We are given the characteristic numbers $e(F)[\mbb{CP}_2]=1$ and $e(F)[\overline{\mbb{CP}}_2]=-1$ by Theorem \ref{shizhouthm} above. 
Using these we compute
$$\< e(F_3^7), \ti g_*[M]\>
=c-d.$$ 
Also
$$\< e(F_3^7), \ti g_*[M]\>=\< \ti g^* e(F), [M]\>=\< e(TM),[M]\>=\chi_M.$$
Equating these two results and eliminating the variable $d=c-\chi_M$ yields the lemma.  
\end{proof}

Note that this is also true with integer coefficients since 
the fourth homology of $G_3^+\mbb R^7$ contains no torsion by \cite{atg2}. 
Next we use this to compute an intersection number.

\beg{lem}\label{intass1} The intersection number ~$\tilde g_*[M]\bullet[ASS]=c.$  \end{lem}

\beg{proof} Making use of the above information we compute the following.
$${
\renewcommand{\arraystretch}{2} 
\begin{array}{rcl} 
\tilde g_*[M]\bullet[ASS] 
& = & (c [\mbb{CP}_2]+(c-\chi_M)[\ol{\mbb{CP}}_2]) \bullet[ASS]\\
& = & c \< [ASS], \tn{PD}^{-1} [\mbb{CP}_2]\> + 
(c-\chi_M)\< [ASS], \tn{PD}^{-1} [\ol{\mbb{CP}}_2] \> \\
& 
= & c\, \< [ASS], {p_1E+eF\over 2}eF \> + 
(c-\chi_M)\, \< [ASS], {p_1E-eF\over 2}eF\>\\
& = & (c-\chi_M/2)\, \< [ASS], p_1EeF \>
+\chi_M/2\, \< [ASS], e^2F \>\\
& 
= & (c-\chi_M/2) \cdot 1 
+\chi_M/2\, \< [ASS], p_1^2E \>\\
& = & (c-\chi_M/2)+\chi_M/2\\
& = & c,
\end{array} }$$  

\noindent where we have used the Theorem \ref{shizhouthm} again. 
We used part $(d)$ by dualizing the projective spaces for illustration purposes but employing $(c)$ and dualizing associatives again would also suffice.  
Here $[ASS]$ denotes the assocative subspace of the Grassmannian. This is the subset of the Grassmannian which correspond to the 3-planes in $\mbb R^7$ calibrated by the 3-form $\phi$. \end{proof}

\noindent Rather than pushing the chain forward, we can use pullbacks of forms as well again to recompute the same number.

\beg{lem}\label{intass2} Alternatively we can compute 
the intersection number as ~$$\tilde g_*[M]\bullet[ASS]={1\over 2}(\chi-3\tau).$$ 
Also we compute 
$$\tilde g_*[M]\bullet[\wt{ASS}]=-{1\over 2}(\chi+3\tau).$$
 \end{lem}

\beg{proof} We use the Theorem \ref{shizhouthm} so that, 
$${
\renewcommand{\arraystretch}{2} 
\begin{array}{rcl} 
\tilde g_*[M]\bullet[ASS] 
& = & \< \,\tn{PD}^{-1}[ASS] , \tilde g_*[M]\>    \\
& = & \< {1\over 2}(p_1E+eF) , \tilde g_*[M]\>    \\
& = & {1\over 2}\, \< p_1(*^*E)+e(*^*F) ,  g_*[M]\>    \\
& = & {1\over 2}\, \< p_1F + eE , g_*[M]\>    \\
& 
  = & {1\over 2}\, \< -p_1E + eE , g_*[M]\>    \\
& = & {1\over 2}\, \< -p_1(TM) + e(TM) , [M]\>.    \\
\end{array} }$$  
because $E\oplus F=\underline{\mbb R}^7$ is trivial, $g^*(E)=TM$ 
and the Hirzebruch signature formula 
\cite{hirzebruchsignaturetheorem}. Reversed oriented associative Grassmannian case is 
similar. 
\end{proof}

\noindent Comparing the two Lemmata we obtain the value $c={1\over 2}(\chi-3\tau)$. 
This improves our Lemma \ref{c1} and gives us the full fundamental class formula as follows. 

\begin{prop}\label{image} 
The image of the canonical class in $H_4(G_3^+\mbb R^7;\mbb R)$ is given by
$$\tilde g_*[M]= {1\over 2}(\chi-3\tau) \, [\mbb{CP}_2]-{1\over 2}(\chi+3\tau)\,  [\overline{\mbb{CP}}_2]$$ for any immersion $i: M\longrightarrow \mbb R^7$.
\end{prop}
\ni In particular, this tells us that the image of the orientation class 
is independent from the immersion. This result is used in the proof of  
a main result in \cite{ibram3}. If $\chi=\tau=0$ then according to our Proposition \ref{image} intersection numbers of the image of the canonical class 
with coassociatiatives is automatically zero. 
Now, from this point on assume that the immersed 4-manifold is coassociative-free. 
We are ready to prove our vanishing result. 

\begin{thm}[Vanishing] \label{vanishingtheorem} If $M^4$ is closed and $i: M\to \mbb R^7$ is a
coassociative($*\phi_0$)-free immersion, then the Euler characteristic $\chi_M$ and 
the signature $\tau_M$ vanishes.
\end{thm}
\beg{proof} Since through the star map we have $$g_*[M]\bullet[COASS]=*_* \circ g_*[M]\bullet*_*[COASS]=\tilde g_*[M]\bullet[ASS],$$ the answer $c$ above is the intersection of the image of the tangent planes with the coassociative planes in $G_4^+\mbb R^7$. 
Free condition implies that the intersection numbers of tangent planes of $M$ with coassociative and 
reversed oriented coassociatives are zero. These are the intersection numbers in 
Lemma \ref{intass2} and gives the linear system
$${
\renewcommand{\arraystretch}{2} 
\begin{array}{rcc}
0&=& {1\over 2}(\chi-3\tau) \\
0&=&-{1\over 2}(\chi+3\tau).\end{array}}$$
\end{proof}

\vspace{.05in}

\ni This is in accordance with the 
following generalization of a theorem of \cite{ibram}.  
This 
Theorem 
tells us that the Euler characteristic 
 is zero even for any manifold with $G_2$ structure instead of the Euclidean 7-space. 

\begin{thm}[$\chi$ vanishing] \label{ibramgenel} Let $i: M \rightarrow X$ be a coassociative-free immersion of a smooth 4-manifold into a 7-manifold $(X,\varphi)$ with a $G_2$ structure. Then the Euler characteristic $\chi_M$ vanishes. 
\end{thm}

\ni The same proof carries on if one takes the nowhere vanishing three form, 
$$\eta:=i^*(\varphi|_{i(M)})\neq 0$$ 
\ni provided by $*\varphi$ freeness. 
Then one have to use an arbitrary metric to take its Hodge star  
$*\eta\in\Lambda^1M$ and convert to a nowhere vanishing vector field on $M$ by the metric duality. 

\section{Cayley-free Embeddings} 
\label{seccayleyfree}

The intersection theoretic computation techniques that we used in the previous section \label{embeddings} can also be used in the Cayley case as well. This case is easier because the fundamental class formula is known already. 
Let $f: M\longrightarrow \mbbr^8$ be an immersion of a compact oriented 
4-manifold, 
then its Gauss map $g: M\longrightarrow G_4^+\mbbr^8$ is computed as 
\cite{shizhou},
$$g_*[M]={1 \over 2}\chi[G(4,5)]+\lambda[G(1,5)]+{3 \over 2}\tau[G(2,4)].$$

\ni Here $\lambda={1\over 2}g^*eF[M]$ and 
$\tau=\tau(M)={1\over 3}g^*p_1E[M]={1\over 3}p_1[TM]$ is the signature. We will be dealing with the embedding case, so that $\lambda=0$.     
We intersect the Gauss image of the embedded 4-manifold with the Cayley and anti-Cayley planes. We can compute the first intersection number as follows,
$${\renewcommand{\arraystretch}{2} 
\begin{array}{rcl} 
g_*[M]\bullet[CAY] 
& = & \< \,\tn{PD}^{-1}[CAY] , g_*[M]\> \\ 
& = & \<{1\over 2}(p_1E-eE+eF) , {\chi\over 2}\,[G(4,5)]+{3\tau\over 2}\,[G(2,4)]\> \\
& = & -{1 \over 2}(\chi-3\tau)\\
\end{array} }$$ 

\ni using the integration table in page 517. To be able to compute the intersection number with the negative Cayley locus, we have to figure out the Poincar\' e dual. 
Reversing the orientation of the Cayley planes, orientation of the bundles are reversed, so that the sign of the Euler classes are changed however the first 
Pontrjagin class does not change sign. So that the Poincar\' e dual becomes, 
$$\tn{PD}^{-1}[\wt{CAY}] = {1\over 2}(p_1E+eE-eF).$$
Inserting this, we compute the second intersection number as follows, 
$${\renewcommand{\arraystretch}{2} 
\begin{array}{rcl} 
g_*[M]\bullet[\wt{CAY}] 
& = & \< \,\tn{PD}^{-1}[\wt{CAY}] , g_*[M]\> \\ 
& = & {1 \over 2}(\chi+3\tau).\\
\end{array} }$$ 

\ni Cayley-free condition implies that these two intersection numbers has to be zero, so the linear equations implies the vanishing similarly as in the previous section. 

\begin{thm}[Vanishing] \label{cayleyfreevanishingtheorem} 
If $M^4$ is compact and $i: M\to \mbb R^8$ is a
Cayley($\psi_0$)-free embedding, then the Euler characteristic $\chi_M$ and 
the signature $\tau_M$ vanishes.
\end{thm}

\ni Using the appropriate h-principle techique, this result can be extended to the 
Cayley free embeddings into $8$-manifolds with $Spin_7$ structure. See \cite{ibram4} for details.

\section{
Contractibility of the Gauss map 
}
\label{hurewiczinjectivity}

In this section we will 
focus on the Gauss map, and prove parallelizability 
through it. We shall start with proving that the Hurewicz homomorphism
$$h_n : \pi_n (G_3^+\mbb R^7)\longrightarrow H_n(G_3^+\mbb R^7;\mbb Z)$$
is injective at the level $n=4$. This homomorphism is defined by sending a homotopy class to its homology class.  
So that $h_n([f])=f_*[S^n]$ where $f$ represents a homotopy class, and $f_*$ is 
the push forward at the homology level.

\vspace{.05in}

The idea of the proof is to use the generalized {\em Mod-$\mc C_p$ Hurewicz Theorems}. 
We start with an introduction to them which follows \cite{daviskirk}. Alternative classical resources are \cite{moshertangora} and \cite{hu}.
For a subset $P$ of prime numbers, let $\mc C_P$ denote the class of 
torsion abelian groups which has no elements of order a positive power of
$p\in P$. This class actually satisfies the properties of a so-called {\em Serre Class}. If we take $P=\{p\}$ then we just use the notation $\mc C_p$.
So as an example we have 
$$\mc C_7=\tn{Torsion 
abelian groups which has no element of order} ~7^k, ~k\in \mbb Z^+.$$
Obviously the groups $\mbb Z_5, \mbb Z_{24}$ etc. but $\mbb Z_{49}$ are in this category. 
A homomorphism $\varphi : A\to B$ between two abelian groups is called
a {\em $\mc C_P$-monomorphism} if $\tn{Ker}\,\varphi\in \mc C_P$,
a {\em $\mc C_P$-epimorphism} if $\tn{Coker}\,\varphi\in \mc C_P$ and 
a {\em $\mc C_P$-isomorphism} if both of these conditions are satisfied. 
Now we are ready to state the following classical theorem whose proof involves spectral sequences \cite{daviskirk}.

\beg{thm}[Mod-$\mc C_P$ Hurewicz]\label{modphurewicz} 
Let $X$ be $1$-connected and 
$\pi_i(X)\in\mc C_P$ for all $i<n$. 
Then, $H_i(X;\mbb Z)\in \mc C_P$ for all $0<i<n$ and the
Hurewicz map $h_n: \pi_n(X)\to H_n(X;\mbb Z)$ is a $\mc C_P$-isomorphism. \end{thm}

As an application we obtain a central result of this section.

\beg{lem}[Injectivity] The Hurewicz homomorphism
$$h_4 : \pi_4 (G_3^+\mbb R^7)\longrightarrow H_4(G_3^+\mbb R^7;\mbb Z)$$
is injective. \end{lem}

\beg{proof} Following \cite{atg2}, 
$G_3^+\mbb R^7$ is simply-connected so $1$-connected, and
possessing $\pi_{0123}=\{0, 0, \mbb Z_2, 0\}$ as the first four homotopy groups
none of which contains elements of order $7^k, ~k\in\mbb Z^+$ hence of class $\mc C_7$. 

Taking $X=G_3^+\mbb R^7$ and $n=4$ 
in the above Mod-$\mc C_7$ Hurewicz Theorem \ref{modphurewicz}, we get that the
Hurewicz map $h_4: \pi_4(X)\to H_4(X;\mbb Z)$ is a $\mc C_7$-isomorphism. 
In particular $\tn{Ker}\,h_4  \in  \mc C_7$.
Besides that $\tn{Ker}\,h_4$ is certainly contained in 
$\pi_4 (G_3^+\mbb R^7)=\mbb Z \oplus \mbb Z$ from \cite{atg2} which is torsion free.
$\mc C_7$ is a class of torsion groups, and the only torsion subgroup
of $\mbb Z \oplus \mbb Z$ is the trivial one, hence $\tn{Ker}\,h_4=0$.
\end{proof}
\ni It is a curious question whether this map is surjective as well. We leave it as an exercise to the reader. Covering all homology classes by means of spheres does not seem likely in a middle level  though. 

\vspace{.05in}

Combining the Proposition \ref{image} and the vanishing Theorem \ref{vanishingtheorem} 
we obtain the following consequence.

\beg{cor}If $M^4$ is compact and $i: M\to \mbb R^7$ is a
coassociative($*\phi_0$)-free immersion, then $g_*[M]=0$ in $H_4(G_4^+\mbb R^7;\mbbz)$ for the associated Gauss map $g$.
\end{cor}


\noindent This observation eventually leads to the contractibility of the Gauss map.

\beg{thm}\label{contractible} If $M^4$ is a compact, oriented, spin manifold 
and $i: M\to \mbb R^7$ is a
coassociative ($*\phi_0$) free immersion, then its image $g(M)$ under the Gauss map $g: M\to G_4^+\mbb R^7$ is contractible, so that $M$ is parallelizable. 
\end{thm}
\beg{proof} The proof uses the obstruction theory. See \cite{hu,milnor}. 
We will work with the equivalent map $\tilde g : M\to G_3^+\mbb R^7$ for convenience. 
We shrink the map skeleton by skeleton. Restriction of $g$ to the 0-th and
1-st skeleton of $M$ is easily contracted to a point by a homotopy because the underlying space $G_4^+\mbb R^7$ is connected and simply connected. Next comes the problem of shrinking the map over the 2-skeleton.
$$\tilde g : M_{(2)} / M_{(1)} \to G_3^+\mbb R^7.$$
Since we already shrinked the map over the 1-skeleton, this gives a cocycle 
hence a cohomology class 
in the second cohomology of $M$ with $\pi_2$ coefficients.
$$\mathfrak{o}_2 \in  H^2( M ; \{ \pi_2(G_3^+\mbb R^7) \} ) = H^2(M;\mbbz_2).$$
The Stiefel-Whitney class $w_2$ is equal to the obstruction $\mathfrak{o}_2$. Since $M$ is spin, by our assumption we get rid of this obstruction and the next one is
$$\mathfrak{o}_3 \in  H^3( M ; \{ \pi_3(G_3^+\mbb R^7) \} ) = 0,$$
since the homotopy group $\pi_3(G_3^+\mbb R^7)$ is trivial as computed in \cite{atg2}. 
The last obstruction lies
$$\mathfrak{o}_4 \in  H^4( M ; \{ \pi_4(G_3^+\mbb R^7) \} ) = 
H^4(M; \{\mbbz\oplus\mbbz\}).$$ 
Since $M$ is a smooth 4-manifold and the top homology is either $\mbbz$ or trivial, we can use a cell complex decomposition with only one 4-cell. 
After we have contracted the 3-skeleton of $M$ to a point in $G_3^+\mbb R^7$, 
the 4-cell gives a map $S^4 \longrightarrow G_3^+\mbb R^7$, homology class of which is the same as $\tilde g_*(M)$. We have already shown 
that the homology class 
$[\tilde g_*(M)]$ is trivial in $H_4(G_3^+\mbb R^7 ;\mbbz)$. In Section 
\ref{hurewiczinjectivity} we have shown that the Hurewicz homomorphism 
$$h_4 : \pi_4 (G_3^+\mbb R^7)\longrightarrow H_4(G_3^+\mbb R^7;\mbb Z)$$
is injective. Hence the homotopy class of the 4-cell map is trivial, consequently $\tilde g(M)$ is contractible. Hence $M$ has trivial tangent  
bundle. 
\end{proof}


\section{Some examples}\label{secexamples}

In this section we will give some examples to illustrate our theorems. 
In particular when we lift some assumptions we will see that there are spaces 
which become no longer embeddable in the appropriate way. 
Our main theorem applies to some connected sums of the panelled web 4-manifolds 
denoted by $M^1_n,M^2_{g,n},M^3_{g,n},M^4_n$ which are 
constructed in \cite{handle} and also $M^5_{g,n}$ is defined  in \cite{sd}. These manifolds are constructed using some special type of Kleinian groups that goes under the same name. They actually come up with a Riemannian metric which is locally conformally flat. But we are only interested in their underlying smooth structure here. Among their topological invariants, their Euler characteristics are computed as follows. 
$$\chi(M^1_n)=-4g, ~~~~\chi(M^2_{g,n})=\chi(M^3_{g,n})=\chi(M^5_{g,n})=4-4g-4n, 
~~~~\chi(M^4_n)=-2n.$$
\ni Since these spaces come up with locally conformally flat metrics their signature is zero. 
\vspace{.05in}
The following examples satisfy all the hypothesis of our Theorem \ref{contractible} and Theorem \ref{ibramgenel}. However their 
tangent bundle is non-trivial because of the non-triviality of the signature obstruction $p_1[M]=3\tau(M)\neq 0$. 
In the following Corollary, $K\!3$ denotes the underlying smooth manifold of a smooth quartic in complex projective space, namely the {\em K\!3 Surface}.

\vspace{.03in}

\beg{cor}\label{corexamplenonzerosignature} The following families of 4-manifolds 
\beg{enumerate}
\item $M^1_{11k} \, \sharp\, 2k K\!3$ ~for all~ $k>0$
\item $ M^{2,3,5}_{g,11k-g} \, \sharp\, 2k K\!3$ ~for all~ $11k>g>0$
\item $ M^4_{11k-2} \, \sharp\, k K\!3$ ~for all~ $k>0$
\end{enumerate}
have the invariants $w_{1234}=0$, $\chi=0$ but are not parallelizable, so they do not admit any coassociative-free immersions into $\mbb R^7$, as well as Cayley-free immersions into $\mbbr^8$.
\end{cor}
\beg{proof}We will check the invariants of the $4$-manifolds. Let us start form the Euler characteristic and signature. Using the connected sum formula for the Euler characteristic we obtain the following. 
$$\chi(M^1_{11k} \, \sharp\, 2k K\!3)=\chi(M^1_{11k})-2k\chi(S^4) +2k\chi(K\!3)
=-4\cdot 11k-2k\cdot 2+2k\cdot 24=0.$$
\ni Orientability and being spin is preserved under the connected sum operation. Since all the building blocks we use are orientable and spin, so their connected sum hence the 
first and second Stiefel-Whitney classes vanish: $w_1=0$ and $w_2=0$ as obstructions. The first Steenrod operator coincides with the Bockstein homomorphism \cite{moshertangora} and applying the Wu's explicit formula \cite{milnor} and orientability we get 
$$\beta w_2= Sq^1w_2=w_1w_2+w_3=w_3.$$
See also \cite{hatchervectorbundles}. Because this is a homomorphism, $w_3=0$ as well. We have already computed the Euler characteristic as zero, so its mod $2$ reduction, hence $w_4=0$. 
The signature is additive under the connected sum operation, so,
$$\tau(M^1_{11k} \, \sharp\, 2k K\!3)=\tau(M^1_{11k})+2k\tau(K\!3)
=0+2k\cdot(-16)=-32k,$$
\ni which is nontrivial, so that the tangent bundle is nontrivial as well. 
\end{proof}
\ni Consequently these examples show that the vanishing of signature is a crucial 
necessary condition for coassociative free and Cayley free immersions/embeddings of $4$-manifolds into manifolds with $G_2$ and $Spin_7$-structure, respectively.  
The panelled web manifolds with even indices have vanishing signature and other invariants. On the other hand, they have strictly negative Euler characteristic, hence they do not embed in a free way to any manifold with $G_2$ or $Spin_7$ structure.  
Other family of examples can be constructed using surgeries, 
like infinite family of homotopy $K\!3$ surfaces, for example knot surgered 
symplectic homotopy $K\!3$ surfaces $E(2)_K$, where $K$ is any fibered knot. 
See \cite{chenkwasiksymplectichomotopyK3surfaces} and references therein for an  
overview and current results in the subject. 

\vspace{.05in}

We also give examples of parallelizable 4-manifolds with 
complicated non-simply-connected topology. 
Here $S^2\times S^2$ stands for the $4$-manifold which is the product of two 2-spheres. 

\beg{thm}The following families of 4-manifolds are parallelizable. 
\beg{enumerate}
\item $M^1_g \, \sharp\, 2g S^2\times S^2$ ~for all~ $g>0$
\item $ M^{2,3,5}_{g,n} \, \sharp\, (2g+2n-3) S^2\times S^2$ ~for all~ $g,n>0$
\item $ M^4_n \, \sharp\, (n-1) S^2\times S^2$ ~for all~ $n>0$.
\end{enumerate}
\end{thm}

\beg{proof}Again we are supposed to check the invariants. Taking into consideration that $\chi(S^2\times S^2)=4$ and $\tau(S^2\times S^2)=0$ the proof is similar to that of the Corollary \ref{corexamplenonzerosignature}. 
One can check their invariants as $w_{1234}=e=p_1=0$. So by the classical parallelizability Theorem \ref{thmparallelizabilityHHM} these $4$-manifolds all have trivial tangent bundle. 
\end{proof}

Consequently they satisfy the necessary conditions of our vanishing Theorems \ref{vanishingtheorem} and \ref{cayleyfreevanishingtheorem}.   
Then combining with the results in \cite{ibram3} and \cite{ibram4} 
we can conclude that they have coassociative-free and Cayley-free embeddings into 
$\mbb R^7$ or $\mbb R^8$ and other 
manifolds with $G_2$ structure 
or $Spin_7$ structure which are flat in a neighborhood of a point, respectively. 
Hence these $4$-manifolds with arbitrarily large fundamental groups (also arbitrarily large first or second Betti number) are freely embeddable into flat tori like $\mbb T^7$ or $\mbb T^8$.


\bigskip


{\small
\beg{flushleft}
\textsc{Orta mh. Z\"ubeyde Han\i m cd. No 5-3 Merkez 74100 Bart\i n
, T\" urk\'{i}ye.}\\
\textit{E-mail address:} \texttt{\textbf{kalafat@\,math.msu.edu}}
\end{flushleft}
}

\newpage

\vspace{1cm}


\bibliography{coassf}{}

\begin{thebibliography}{KOA13}

\bibitem[AK12]{handle}
Selman Akbulut and Mustafa Kalafat.
\newblock A class of locally conformally flat 4-manifolds.
\newblock {\em New York J. Math.}, 18:733--763, 2012.

\bibitem[AK16]{atg2}
Selman Akbulut and Mustafa Kalafat.
\newblock Algebraic topology of {$G_2$} manifolds.
\newblock {\em Expo. Math.}, 34(1):106--129, 2016.

\bibitem[CK11]{chenkwasiksymplectichomotopyK3surfaces}
Weimin Chen and Slawomir Kwasik.
\newblock Symmetric symplectic homotopy {$K3$} surfaces.
\newblock {\em J. Topol.}, 4(2):406--430, 2011.

\bibitem[DK01]{daviskirk}
James~F. Davis and Paul Kirk.
\newblock {\em Lecture notes in algebraic topology}, volume~35 of {\em Graduate
  Studies in Mathematics}.
\newblock American Mathematical Society, Providence, RI, 2001.

\bibitem[Fan02]{fang}
Fuquan Fang.
\newblock Orientable 4-manifolds topologically embed into {${\Bbb R}^7$}.
\newblock {\em Topology}, 41(5):927--930, 2002.

\bibitem[Hat18]{hatchervectorbundles}
Allen Hatcher.
\newblock {\em Vector Bundles and K-Theory}.
\newblock Book available online at ~
  \texttt{http://www.math.cornell.edu/$\sim$hatcher/}, 2018.

\bibitem[HH58]{hirzebruchhopf}
Friedrich Hirzebruch and Heinz Hopf.
\newblock Felder von {F}l\"achenelementen in 4-dimensionalen
  {M}annigfaltigkeiten.
\newblock {\em Math. Ann.}, 136:156--172, 1958.

\bibitem[Hir95]{hirzebruchsignaturetheorem}
Friedrich Hirzebruch.
\newblock {\em Topological methods in algebraic geometry}.
\newblock Classics in Mathematics. Springer-Verlag, Berlin, 1995.
\newblock Translated from the German and Appendix One by R. L. E.
  Schwarzenberger, With a preface to the third English edition by the author
  and Schwarzenberger, Appendix Two by A. Borel, Reprint of the 1978 edition.

\bibitem[HL82]{harveylawson}
Reese Harvey and H.~Blaine Lawson, Jr.
\newblock Calibrated geometries.
\newblock {\em Acta Math.}, 148:47--157, 1982.

\bibitem[HL09]{harveylawsonpotential}
F.~Reese Harvey and H.~Blaine Lawson, Jr.
\newblock An introduction to potential theory in calibrated geometry.
\newblock {\em Amer. J. Math.}, 131(4):893--944, 2009.

\bibitem[Hu59]{hu}
Sze-tsen Hu.
\newblock {\em Homotopy theory}.
\newblock Pure and Applied Mathematics, Vol. VIII. Academic Press, New
  York-London, 1959.

\bibitem[Joy00]{joyce}
Dominic~D. Joyce.
\newblock {\em Compact manifolds with special holonomy}.
\newblock Oxford Mathematical Monographs. Oxford University Press, Oxford,
  2000.

\bibitem[KOA13]{sd}
Mustafa Kalafat, Y{\i}ld{\i}ray Ozan, and H\"ulya Arg\"uz.
\newblock Self-dual metrics on non-simply connected 4-manifolds.
\newblock {\em J. Geom. Phys.}, 64:79--82, 2013.

\bibitem[Mas58]{masseyparallelizability}
W.~S. Massey.
\newblock On the cohomology ring of a sphere bundle.
\newblock {\em J. Math. Mech.}, 7:265--289, 1958.

\bibitem[MS74]{milnor}
John~W. Milnor and James~D. Stasheff.
\newblock {\em Characteristic classes}.
\newblock Princeton University Press, Princeton, N. J.; University of Tokyo
  Press, Tokyo, 1974.
\newblock Annals of Mathematics Studies, No. 76.

\bibitem[MT68]{moshertangora}
Robert~E. Mosher and Martin~C. Tangora.
\newblock {\em Cohomology operations and applications in homotopy theory}.
\newblock Harper \& Row, Publishers, New York-London, 1968.

\bibitem[SZ14]{shizhou}
Jin Shi and Jianwei Zhou.
\newblock Characteristic classes on {G}rassmannians.
\newblock {\em Turkish J. Math.}, 38(3):492--523, 2014.

\bibitem[Tho68]{emeryparallelizability}
Emery Thomas.
\newblock Vector fields on low dimensional manifolds.
\newblock {\em Math. Z.}, 103:85--93, 1968.

\bibitem[\"U11]{ibram}
\.Ibrahim \"Unal.
\newblock Topology of {$\phi$}-convex domains in calibrated manifolds.
\newblock {\em Bull. Braz. Math. Soc. (N.S.)}, 42(2):259--275, 2011.

\bibitem[\"U15]{ibram3}
\.Ibrahim \"Unal.
\newblock {$h$}-principle and {$\phi$}-free embeddings in calibrated manifolds.
\newblock {\em Internat. J. Math.}, 26(7):1550052, 16, 2015.

\bibitem[\"U18]{ibram4}
\.Ibrahim \"Unal.
\newblock A note on the {G}auss maps of {C}ayley-free embeddings into
  {S}pin(7)-manifolds.
\newblock {\em Differential Geom. Appl.}, 61:1--8, 2018.

\end{thebibliography}
\bibliographystyle{alphaurl}
\end{document}